%SECOND BATCH OF HERTLING'S CORRECTIONS INSERTED
%PROOF FROM [BAMA] INSERTED
\input amstex
\documentstyle{amsppt}
\document
\magnification=1200
\NoBlackBoxes
\vsize17cm
\hoffset=0.7in
\voffset=1in
\nologo

%\hfill{\it file Amywork/updateF.tex, version 19.03.08}

\bigskip

\centerline{\bf AN UPDATE }

\medskip

\centerline{\bf ON SEMISIMPLE QUANTUM COHOMOLOGY}

\medskip

\centerline{\bf AND $F$--MANIFOLDS}

\medskip

\centerline{\bf Claus Hertling, Yuri I. Manin, Constantin Teleman}

\medskip

%\centerline{\it Max--Planck--Institut f\"ur Mathematik, Bonn, Germany,}

%\centerline{\it and Northwestern University, Evanston, USA}

\bigskip

{\bf Abstract.}  In the first section of this note we show that the Theorem 1.8.1 of
Bayer--Manin ([BaMa]) can be strengthened in the following way: {\it if
the even quantum cohomology of a projective algebraic manifold $V$ is generically
semi--simple, then $V$ has no odd cohomology and is of Hodge--Tate type.}
In particular, this addressess a question in [Ci].

In the second section, we prove that {\it an analytic (or formal) 
supermanifold $M$ with a given supercommutative
associative $\Cal{O}_M$--bilinear multiplication on its tangent sheaf
$\Cal{T}_M$ is an $F$--manifold in the sense of [HeMa], iff its spectral
cover as an analytic subspace of the cotangent bundle $T^*_M$
is coisotropic of maximal dimension.} This answers a question of V.~Ginzburg.

Finally, we discuss these results in the context of mirror
symmetry and Landau--Ginzburg models for Fano varieties.

\bigskip

\centerline{\bf \S 0. Introduction}

\medskip

{\bf 0.1. Contents of the paper.} Semisimple Frobenius manifolds have many nice properties: see e.~g. [Du], [Ma], [Te],
[Go1], [Go2], and references therein.
It is important to understand as precisely as possible, which
projective algebraic manifolds $V$ have (generically) semi--simple quantum cohomology.
In this case the quantum cohomology is determined
by a finite amount of numbers, and a mirror 
(Landau--Ginzburg model) can in many cases
be described  explicitly.

\smallskip

If $V$ has non--trivial odd cohomology, its full quantum cohomology
cannot be semi--simple, but its even part is a closed
Frobenius subspace, and in principle it can be semisimple.
In [BaMa], Theorem 1.8.1, it was proved that
if $H^{ev}_{quant}(V)$ is generically semisimple, then
$h^{p,q}(V)=0$ for $p+ q \equiv 0\,\roman{mod}\,2$, $p\ne q.$
In the first section of this note we show that in this case  $h^{p,q}(V)=0$
for  $p+ q \equiv 1\,\roman{mod}\,2$ as well.

\smallskip

Thus, the Theorem 1.8.1 of
Bayer--Manin ([BaMa]) can be strengthened in the following way: {\it if
the even quantum cohomology of a projective algebraic manifold $V$ is generically
semi--simple, then $V$ has no odd cohomology and is of Hodge--Tate type.}
In particular, for the 47 families of Fano threefolds with $b_3(V)>0$, pure
even quantum cohomology cannot be semi--simple, cf. [Ci], p. 826. 

\smallskip

The second section is dedicated to a previously
unpublished result of C.~Hertling (letter dated March 09, 2005,
where it was stated for the pure even case). It shows that an analytic (or formal) 
supermanifold $M$ with a given supercommutative
associative $\Cal{O}_M$--bilinear multiplication on its tangent sheaf
$\Cal{T}_M$ is an $F$--(super)manifold in the sense of [HeMa], iff its spectral
cover as an analytic subspace of the cotangent bundle $T^*_M$
is coisotropic of maximal dimension.

\smallskip

 This answered a question 
posed to Yu.~Manin by  V.~Ginzburg.

\medskip

{\it Acknowledgement.} We are grateful to Arend Bayer for 
illuminating comments on the Proposition 1.2 and for sharing
with us his version of  Dubrovin's
conjecture.

\bigskip

\centerline{\bf \S 1. Semisimple quantum cohomology}

\smallskip

\centerline{\bf and Dubrovin's conjecture}

\medskip

{\bf 1.1. Notation.} Let $V$ be a projective manifold over $\bold{C}$.
We denote by $H^{ev}_{quant}(V)$ its even quantum cohomology
ring. As in [BaMa] and [Ba], it is a topological commutative algebra.
 Multiplication
in it  (the classical cup product plus ``quantum corections'')  is denoted $\circ$. 
The space $H^{ev}(V)$ is embedded
in it as a field of flat vector fields on the respective formal Frobenius manifold.

\smallskip

{\bf 1.2.  Proposition}. {\it If $H^{ev}_{quant}(V)$ is generically reduced,
i.e. has no nilpotents at (the local ring of) the generic
point, then
$H^{odd}(V)=0$.}

\smallskip

{\bf Proof.} Assume that $H^{odd}(V)\ne 0$. Let $\Delta$ be a non--zero class of an odd dimension.

\smallskip

First, we have $\Delta\circ \Delta=0$. In fact,
$\Delta\cup \Delta=0$, because
the cup produt is supercommutative. The quantum
corrections vanish, because the correlators $\langle...\rangle$
are also supercommutative in their arguments, so $\langle\Delta\Delta\Delta^{\prime}...\rangle
=0$. This follows from the fact
that the quantum correlators come from the
$S_n$--covariant maps $H^*(V)^{\otimes n}\to H^*(\overline{M}_{0,n})$,
induced by algebraic correspondences (push--forwards of
virtual fundamental classes).
Covariance holds with respect to the action of $S_n$ on the tensor power
permuting factors and introducing signs as usual
in $Z_2$-graded setting. On the target it renumbers
points and hence leaves the fundamental class invariant.

\smallskip

Now, find another
(odd) class $\Delta^{\prime}$ such that $g(\Delta ,\Delta^{\prime})=1$, where $g$ is the Poincare
form. Then we have $1=g(\Delta,\Delta^{\prime})=g(\Delta\circ \Delta^{\prime}, e)$ where $e$ is the identity
in quantum cohomology. Hence $\Delta\circ \Delta^{\prime}\in
H^{ev}_{quant}(V)$ must be generically non-zero.
But its square is zero because of the first remark.
This contradicts the generic absence of nilpotents in $H^{ev}_{quant}(V)$.
\medskip

{\bf 1.3. Theorem.} {\it  If
the even quantum cohomology of a projective algebraic mani\-fold $V$ is generically
semi--simple, then $V$ has no odd cohomology and is of Hodge--Tate type.}

\smallskip

{\bf Proof.} From the generic semisimplicity and the Proposition 1.2
it follows, that $h^{p,q}(V)=0$ for $p+q\equiv 1\,\roman{mod}\,2.$

\smallskip

To prove that $h^{pq}(V)=0$ for $p+q \equiv 0\,\roman{mod}\,2$ and $p\ne q$,
we reproduce a short reasoning from [BaMa]. It compares the Lie
algebra of Euler vector fields in the semi--simple case
and in the quantum cohomology case.

\smallskip

Firstly, in he semisimple case each Euler vector field
must be of the form $E=d_0\sum_iu_ie_i+\sum_j c_je_j$, where $d_0$
is a constant (weight of $E$, cf. [Ma1], [Ma2]), and $(u_i)$ are (local)
Dubrovin's canonical coordinates, that is, $e_i:=\partial /\partial u_i$
form a complete system of pairwise orthogonal idempotents
in $H^*_{quant}(V)$. Moreover, $(c_j)$ are arbitrary constants.

\smallskip

From this explicit description it follows directly, that if two
Euler fields of non--zero weights commute, they are
proportional.

\smallskip

On the other hand, if $h^{p,q}(V)\ne 0$ for some 
 $p+q \equiv 0\,\roman{mod}\,2$ and $p\ne q$,
 then $H^*_{quant}(V)$ admits two commuting and
 non--proportional Euler vector fields $E_1,E_2$ of weight 1. Namely,
 in the bihomogeneous (with respect to the $(p,q)$--grading) basis
 of flat vector fields $\Delta_a\in H^{p_a,q_a}(V)$, we can take
 $$
 E_1:=\sum_a (1-p_a)x_a\Delta_a +\sum _{p_b=q_b=1} r_b\Delta_b,
 $$
$$
 E_2:=\sum_a (1-q_a)x_a\Delta_a +\sum _{p_b=q_b=1} r_b\Delta_b.
 $$
  Here $(x_a)$ are dual flat coordinates, and  $-K_V=c_1(\Cal{T}_V)=\sum_b r_b\Delta_b.$

\smallskip

  This completes the proof.
\medskip

{\bf 1.4. Dubrovin's conjecture and related insights.} In [Du] (p. 321) the problem
of characterization of varieties $V$ with semisimple quantum cohomology was formulated
explicitly. It was also stated there that  a necessary condition for such
$V$ is to be Fano. This was disproved by A.~Bayer [Ba], who established
that blowing up points on such  a variety does not destroy semisimplicity.
In particular, not only del Pezzo surfaces have semisimple quantum cohomology,
but arbitrary blowups of $\bold{P}^2$ as well.

\smallskip

A.~Bayer has later conjectured that the maximal length of a semi--orthogonal decomposition of 
$D^b(V)$ must coincide with the generic number
of idempotents in $H^*_{quant}(V)$.

\smallskip

Combining the results of [Ba], of this note, and the  further part of Dubrovin's
conjecture stated on p. 322 of [Du] (cf. also [Z]), one can now guess that a necessary and sufficient
condition for semisimplicity is that $V$ is of Hodge--Tate type, 
whose bounded derived coherent category admits a full exceptional
collection $(E_i)$. Moreover, after adjusting some arbitrary
choices, in this case one should be able to identify the Stokes matrix
of its second structure connection with the matrix $(\chi (E_i,E_j))$.

\smallskip

This last statement is now checked, in particular, for three--dimensional
Fano varieties with minimal cohomology in [Go2]. The reader can find  
there more details and explanations about the involvement 
of the vanishing cycles
in the mirror Landau--Ginzburg model. 

\smallskip

All these constructions reflect some facets of Kontsevich's homological
mirror symmetry program. However, one should keep in mind that in this
note we are concerned almost exclusively with a multiplication
on the tangent bundle, i.e. with the structure of an $F$--manifold
(see below).  In order to  invoke mirror symmetry,
we need also to take in consideration a compatible flat metric.
In quantum cohomology, it comes ``for free'' at the start;
it is multiplication that requires a special construction.
In various contexts relevant for mirror symmetry, the metric 
can be described implicitly by at least 
five different kinds of data which we list here for 
reader's convenience.

\smallskip

(a) Values of the diagonal coefficients of the flat metric
$\sum_i\eta_i (du_i)^2$ in canonical coordinates and values of
their first derivatives $\eta_{ij}$ at a tame semi--simple point.
This is initial data for the
second structure connection (cf. [Ma1], II.3).

\smallskip

(b) Monodromy data for the first structure connection
and oscillating integrals for the deformed flat coordinates
(cf. [Gi], [Du], [Sa] and the references therein).

\smallskip
(c) Choice of one of K.~Saito's primitive forms.

\smallskip
(d) Choice of a filtration on the
cohomology space of the Milnor fiber (M.~Saito, cf. [He2]
and the references therein).

\smallskip
(e) Use of the semi--infinite Hodge structure.
This is a refinement of (c), described by S.~Barannikov ([Bar1], [Bar2]).

\bigskip

\centerline{\bf \S 2. $F$--geometry and  symplectic geometry}

\medskip

{\bf 2.1. $F$--structure and Poisson structure.}
Manifolds $M$ considered in this section
can be $C^{\infty}$, analytic, or formal, eventually with
even and odd coordintes (supermanifolds). The ground field
$K$ of characteristic zero is most often $\bold{C}$ or $\bold{R}$.
Each such manifold, by definition,  is endowed with the structure sheaf $\Cal{O}_M$
which is a sheaf of (super)commutative $K$--algebras,
and the tangent sheaf $\Cal{T}_M$ which is a locally
free $\Cal{O}_M$--module of (super)rank equal to the
(super)dimension of $M$. $\Cal{T}_M$ acts on $\Cal{O}_M$
by derivations, and is  a sheaf of Lie (super)algebras
with an intrinsically defined Lie bracket $[\, ,]$. 

\smallskip

There is a classical notion of {\it Poisson structure} on $M$
which endows $\Cal{O_M}$ as well  with a
Lie bracket $\{\, ,\}$ constrained by a well known  identity.

\smallskip

Similarly, {\it an $F$--structure} on $M$ endows $\Cal{T}_M$
with an extra operation: (super)com\-mutative and associative
$\Cal{O}_M$--bilinear multiplication. We denote it
always  $\circ$ and assume that it is endowed with identity: an even
vector field $e$. Then $\Cal{O}_M$ is embedded in $\Cal{T}_M$
as a subalgebra: $f\mapsto fe.$

\smallskip

Given such a multiplication on the tangent sheaf,
we can define its {\it spectral cover $\widetilde{M}$} which is a closed
ringed (super)subspace (generally not a submanifold) in the cotangent 
(super)manifold 
$T^*M$. In the Grothendieck language, it is
simply the relative affine  spectrum  of the sheaf of algebras
$(\Cal{T}_M,\circ )$ on $M$.

\smallskip

More precisely, consider $Symm_{\Cal{O}_M}(\Cal{T}_M)$ as the sheaf
of algebras of those functions on  the cotangent (super)space $T^*_M$
that are polynomial along the fibres of the projection $T^*_M\to M$.
The multiplication in this sheaf will be denoted $\cdot$.
For example, for two local vector fields $X,Y\in \Cal{T}_M (U)$,
$X\cdot Y$ denotes their product as an element of
$Symm^2_{\Cal{O}_M}(\Cal{T}_M)$. 

\smallskip
Consider the
canonical surjective morphism of  sheaves of $\Cal{O}_M$--algebras
$$
Symm_{\Cal{O}_M}(\Cal{T}_M)\to (\Cal{T}_M,\circ )
$$
sending, say, $X\cdot Y$ to $X\circ Y$. Denote its kernel
by $J(M,\circ )$, and let $\widetilde{M}$ be defined by the sheaf of ideals 
$J(M,\circ )$. 
\smallskip

The spectral cover $\widetilde{M}\to M$ is flat, because
$\Cal{T}_M$ is locally free.

\smallskip

Now we will  describe the structure identities  imposed onto
$\{\,,\}$ on $\Cal{O}_M$, resp. $\circ$ on $\Cal{T}_M$.
To this end, recall the notion of the Poisson tensor. Let generally
$A$ be a $K$--linear superspace (or a sheaf of superspaces) endowed with a $K$--bilinear
multiplication and a $K$--bilinear Lie bracket $[\, ,]$. 
Then for any $a,b,c\in A$ put
$$
P_a(b,c):=[a,bc]-[a,b]c-(-1)^{ab}b[a,c] .
\eqno(2.1)
$$
(From here on, $(-1)^{ab}$ and similar notation refers
to the sign occuring in superalgebra when the two neighboring elements get
permuted.)

\smallskip

This tensor will be written for $A=(\Cal{O}_M, \cdot, \{\,,\})$
in case of the Poisson structure,  and for $A=(\Cal{T}_M, \circ , [\,,])$
in case of an $F$--structure.

\smallskip

We will now present parallel lists of basic properties of 
Poisson, resp. $F$--manifolds.

\medskip

{\bf 2.2. Poisson (super)manifolds.} (i)${}_P$. {\it Structure identity:}
for all local functions $f,g,h$ on $M$
$$
P_f(g,h)\equiv 0 .
\eqno(2.2)
$$

(ii)${}_P$. {\it Each  local
function $f$ on $M$ becomes a local
vector field $X_f$ (of the same parity as $f$) on $M$ via
$X_f(g):=\{f,g\}$.}

 This is a reformulation of (2.2).

\smallskip

(iii)${}_P$. {\it Maximally nondegenerate case: symplectic structure.}
There exist local canonical coordinates $(q_i,p_i)$ such that
for any $f,g$
$$
\{f,g\}=\sum_{i=1}^n (\partial_{q_i}f\partial_{p_i}g
- \partial_{q_i}g\partial_{p_i}f).
$$
Thus, locally all symplectic manifolds of the same dimension are isomorphic.
The local group of symplectomorphisms is, however,
infinite dimensional.

\medskip

{\bf 2.3. $F$--manifolds.} (i)${}_F$. {\it Structure identity:}
for all local vector fields $X,Y,Z,U$
$$
P_{X\circ Y}(Z,U)-X\circ P_Y(Z,U)- (-1)^{XY}Y\circ P_X(Z,U)=0.
\eqno(2.3)
$$

(ii)${}_F$. {\it Each  local
vector field on $M$ becomes a local
function  on the spectral cover $\widetilde{M}$ of $M$}. 
\smallskip

As we already mentioned, generally $\widetilde{M}$ is not  a (super)manifold.
In the pure even case this often happens because of nilpotents in 
$\Cal{O}_{\widetilde{M}}$ and/or singularities. 
In the presence of odd coordinates on $M$ nilpotents
by themselves are always present, but typically they cannot
form an exterior algebra over functions of even coordinates
because ranks do not match.

\smallskip

 A theorem due to Hertling
describes certain important cases when $\widetilde{M}$
is a manifold.

\smallskip

(iii)${}_F$. {\it Maximally nondegenerate case: semisimple $F$--manifolds.}
$\widetilde{M}$  will be a manifold and even an unramified covering
of $M$ in the appropriate ``maximally nondegenerate case'',
namely, when $M$ is pure even, and locally $(\Cal{T}_M,\circ )$ is isomorphic 
to $(\Cal{O}_M^d)$ as algebra, $d=\roman{dim}\,M.$

\smallskip

In this case there exist local canonical coordinates $(u_a)$ (Dubrovin's coordinates)
such that the respective vector fields $\partial_a:=\partial/\partial_a$
are orthogonal idempotents:
$$
\partial_a\circ \partial_a =\delta_{ab}\partial_a .
$$
Thus, locally all semisimple $F$--manifolds of the same dimension are isomorphic.
Local automorphisms of an $F$--semisimple structure
are generated by renumberings and shifts of canonical coordinates:
$$
u_a\mapsto u_{\sigma (a)}+c_a 
$$
so that this structure is more rigid than the symplectic one.

\medskip

{\bf 2.4. Spectral cover as a subspace in symplectic supermanifold.}  
There is a structure of  sheaf of Lie algebras 
on $Symm_{\Cal{O}_M}(\Cal{T}_M)$. It is given by
the Poisson brackets  $\{\,, \}$ with respect to the
canonical (super)symplectic structure on  $T^*_M$.

\smallskip

It is easy to check that the  ideal $J=J (M,\circ )\subset Symm_{\Cal{O}_M}(\Cal{T}_M)$
defining $\widetilde{M}$  in 
this sheaf of supercommutative algebras is generated by
all expressions:
$$
e-1,\quad X\circ Y-X\cdot Y, \quad X,Y \in \Cal{T}_M.
\eqno(2.4)
$$
 \medskip
 
 {\bf 2.5. Theorem.} {\it The multiplication $\circ$
 satisfies the structure identity of $F$--manifolds (2.3),
 iff the ideal $J(M,\circ )$ is stable with respect to
 the Poisson brackets. }
 
  \medskip
 
 {\bf Proof.} From (2.2), one easily infers that stability of an ideal
 in a Poisson algebra with respect to the brackets can be checked on
 any system of generators of this ideal. In our case we choose (2.4).

\smallskip 

Clearly, $\{e-1,e-1\} =0$.

\smallskip

If $X,Y$ are local vector fields, then $\{X,Y\}=[X,Y]$

\smallskip

 We will establish by a  direct computation that for
 all $X,Y,Z,W$ as above,
 $$
 \{X\circ Y-X\cdot Y,\, Z\circ W-Z\cdot W\} \equiv
 $$
 $$
 P_{X\circ Y}(Z,W)-X\circ P_Y(Z,W)-(-1)^{XY}Y\circ P_X(Z,W)\ \roman{mod}\, 
 J(M,\circ )
\eqno(2.5)
$$
and
$$
\{e-1, X\circ Y-X\cdot Y\} = [e,X\circ Y]-X\cdot [e,Y] -[e,X]\cdot Y.
\eqno(2.6)
$$
\smallskip

Assume that this is done. From (2.5) and (2.6) it follows that if (2.3) holds, then
$J(M,\circ )$ is stable with respect to the Poisson brackets.  For (2.6), one uses the identity 
$[e, X\circ Y]=X\circ [e,Y] + [e,X]\circ Y$ which follows from (2.3) by choosing
$X=Y=e$ and renaming $Z,U$.

\smallskip
Conversely, if $J(M,\circ )$ is stable with respect to the brackets, then the right--hand side of (2.5)
must  belong to  $J(M,\circ )$. But it lies
in the degree 1 part of the symmetric algebra of $\Cal{T}_M$,
which projects onto $\Cal{T}_M$. Hence it must vanish,
and as a result, the right hand side of (2.6) must belong
to $J(M,\circ )$ as well.

\smallskip

It remains to check (2.5) and (2.6). We will briefly indicate how to do it,
restricting ourselves to the clumsier case (2.5).

\smallskip

First of all, the right hand side of (2.5) can be rewritten
as follows: 
$$
P_{X\circ Y}(Z,W)-X\circ P_Y(Z,W)-(-1)^{XY}Y\circ P_X(Z,W)=
$$
$$
[X\circ Y,Z\circ W]-[X\circ Y,Z]\circ W-
(-1)^{(X+Y)Z}Z\circ [X\circ Y,W]
$$
$$
-X\circ [Y,Z\circ W] -(-1)^{XY}Y\circ [X,Z\circ W]+ X\circ [Y,Z]\circ W
+(-1)^{YZ} X\circ Z\circ [Y,W]
$$
$$
 + (-1)^{XY}Y\circ [X,Z]\circ W
+(-1)^{X(Y+Z)}Y\circ Z\circ [X,W] .
\eqno(2.7)
$$
It turns out that (2.7) is in fact a tensor, that is $\Cal{O}_M$--polylinear
in $X,Y,Z,W$. See [Me1], [Me2] for a discussion and  operadic generalizations  
of the condition of its vanishing. 
\smallskip

In our context, this formula is convenient, because  a straightforward decomposition of 
the left hand side of (2.5) into Poisson monomials
(constructed using two operations) gives exactly the same
list of monomials as in (2.7) modulo $J(M,\circ )$, with the same signs.

\smallskip

Here are samples of calculations.  

\smallskip

The first term
$\{X\circ Y, Z\circ W\}$ at the left hand side of (2.5) coincides with the
first term in (2.7).

\smallskip

Using the Poisson identity (2.2), we find further:
$$
-\{X\circ Y, Z\cdot W\} =-\{X\circ Y, Z\}\cdot W-
(-1)^{(X+Y)Z}Z\cdot \{X\circ Y, W\}.
$$
Modulo $J(M, \circ )$, this can be replaced by
$$
-[X\circ Y, Z]\circ W-
(-1)^{(X+Y)Z}Z\circ [X\circ Y, W]
$$
which corresponds to the second and third terms
of (2.7).

\smallskip

We leave the rest as an exercise to the reader.

\medskip

{\bf 2.5.1. Reduced spectral cover.} Contrary to what 
might be expected, the condition 
$$
\{J(M,\circ ),J(M,\circ )\}\subset J(M,\circ )
$$
{\it does not} imply the respective condition
for the radical of $J(M,\circ )$ even in the pure even case.
This means that $\widetilde{M}_{red}$ need not
be a Lagrange subvariety, even if it comes from
an $F$--manifold.

\smallskip

This can be shown on the following explicit examples.

\smallskip
We will construct two families of everywhere indecomposable (see 2.6 below) $F$-manifolds
in terms of the ideals $J$, defining (nonreduced)
subspaces $\widetilde{M} \subset T^*M$. In order to give rise to $F$-manifolds
with $\pi :\,\widetilde M\to M$ as their spectral cover, they have to satisfy
the following conditions:

\smallskip

(a) The projection $\widetilde{M}\to M$ is flat of degree $n=\,\roman{dim}\, M$ and the 
canonical map $\Cal{T}_M\to \pi_*(\Cal{O}_{\widetilde{M}})$ is an isomorphism.

\smallskip
To check this by direct calculations, 
we will choose (pure even) local coordinates $(t_1,\dots ,t_n)$ on $M$ in such a way that
$e=\partial /\partial t_1.$ By $(y_1,\dots ,y_n)$ we will denote
the conjugate coordinates along the fibres of $T^*(M)$.

\smallskip

(b) $\{J,J\}\subset J$.

\smallskip

We will see that in these examples
$$
\{\sqrt{J},\sqrt{J}\}\not\subset \sqrt{J}.
$$

\medskip

{\bf 2.5.2. The first family.} Here we put
$$
J=(y_1-1,\; (y_i-\rho_i)(y_j-\rho_j)),
$$
with $\rho_1=1$ and $\rho_i\in \Cal{O}_M$ for $i\geq 2$ such that
$\partial_1 \rho_i=0$.
Clearly, (a) and (b) are satisfied. The radical of $J$ is
$$
\sqrt{J} = (y_1-1; y_2-\rho_2, \dots ,y_n-\rho_n).
$$
We have $\{\sqrt{J},\sqrt{J}\}\not\subset \sqrt{J}$, if
$$
\partial_i \rho_j \neq \partial_j\rho_i
\quad \hbox{ for some }i,j\geq 2 \hbox{ with }i\neq j.
$$
The algebra $T_tM$ at any point $t\in M$ is isomorphic to
$\hbox{\bf C}[x_1,...,x_{n-1}]/(x_ix_j)$.

\medskip

{\bf 2.5.3. The second family.}  Here we put for any $n\ge 3$
$$
J=(y_1-1,\; (y_2-\rho_2)^2,\; (y_2-\rho_2)\cdot y_3, \;
y_3^{n-1},\; y_4-y_3^2,\; y_5-y_3^3,...,\; y_n-y_3^{n-2}),
$$
with
$$
\rho_2(y,t)= t_3y_1+\sum_{k=3}^{n-1} (k-1)t_{k+1}\cdot y_k.
$$
Now, (a) is rather obvious, but checking (b) requires a calculation
which we omit.
The radical of $J$ is
$$
\sqrt{J} = (y_1-1,\; y_2-t_3\cdot y_1,\; y_3,y_4,y_5,...,y_n).
$$

\smallskip
The algebra $T_tM$ at any point $t\in M$ is isomorphic to
$\hbox{\bf C}[x_2,x_3]/(x_2^2,x_2x_3,x_3^{n-1})$.

\smallskip

We will now explain in which context the considerations
of this section can be related to the problems, arising
in the study of semisimple quantum cohomology

\medskip

{\bf 2.6. Hertling's local decomposition theorem.} For any point
$x$ of a pure even $F$--manifold $M$, the tangent space $T_xM$
is endowed with the structure of a $K$--algebra.
This $K$--algebra can be represented as a direct sum of local
$K$--algebras. The decomposition is unique in the following sense:
the set of pairwise orthogonal idempotent tangent
vectors determining it is well defined.

\smallskip

C. Hertling has shown that this decomposition extends to a
neighborhood of $x$. 

\smallskip

More precisely, define the sum of two $F$--manifolds:
$$
(M_1,\circ_1,e_1) \oplus (M_2,\circ_2,e_2) :=
(M_1\times M_2,\circ_1\boxplus\circ_2,e_1\boxplus e_2)
$$
A manifold is called {\it indecomposable} if it cannot be represented
as a sum in a nontrivial way. 

\medskip

{\bf 2.6.1. Theorem.} {\it Every germ $(M,x)$ of a complex
analytic $F$--manifold 
decomposes into a direct sum
of indecomposable germs such that for each summand,
the tangent algebra at $x$ is a local algebra.

\smallskip

This decomposition is unique
in the following sense: the set of pairwise orthogonal idempotent
vector fields determining it is well defined.} 

\medskip

For a proof, see [He], Theorem 2.11. 

\smallskip

Furthermore, we have ([He], Theorems 5.3 and 5.6):

\medskip

{\bf 2.7. Theorem.} {\it (i) The spectral cover space $\widetilde{M}$
of the $F$--structure on the germ of the unfolding space
of an isolated hypersurface singularity is smooth.

\smallskip

(ii) Conversely, let $M$ be an irreducible germ of a generically
semisimple $F$--manifold with the smooth spectral cover $\widetilde{M}.$
Then it is (isomorphic to) the germ of the unfolding space
of an isolated hypersurface singularity. Moreover, any isomorphism
of germs of such unfolding spaces compatible with their
$F$--structure comes from a stable right equivalence of the
germs of the respective singularities.}

\smallskip

Recall that the stable right equivalence is generated by adding sums
of squares of coordinates and making invertible analytic
coordinate changes.

\smallskip

In view of this result, it would be important to understand
the following

\medskip

{\bf 2.8. Problem.} {\it Characterize those varieties $V$ for which the
 quantum cohomology Frobenius spaces $H^*_{quant}(V)$ have
smooth spectral covers.}

\smallskip

Theorem 2.7 produces for such manifolds a weak version
of Landau--Ginzburg model, and thus gives a partial
solution of the mirror problem for them.

\newpage

\centerline{\bf References}

\medskip

[AuKaOr] D.~Auroux, L.~Katzarkov, D.~Orlov. {\it Mirror symmetry for del Pezzo surfaces: vanishing cycles and coherent sheaves.}  Invent. Math.  166  (2006),  no. 3, 537--582.  e-arxiv 0506166

\smallskip

[Bar1] S.~Barannikov. {\it Semi--infinite Hodge structures and mirror symmetry for projective spaces.}
e-arxiv  0010157

[Bar2] S.~Barannikov. {\it Semi--infinite variations of Hodge structure and integrable
hierarchies of KdV type.} Int. Math. Research Notices, 19 (2002), 973--990. e-arxiv math.AG/00148.

\smallskip

[Ba] A.~Bayer.  {\it Semisimple quantum cohomology and blowups.}
Int. Math. Research Notices, 40 (2004), 2069--2083.
 e-arxiv math.AG/0403260.

\smallskip

[BaMa] A.~Bayer, Yu.~Manin. {\it (Semi)simple exercises in 
quantum cohomology.} Proc. of the Fano conference,
Univ. Torino, Turin, 2004, 143--173. e-arxiv math.AG/0103164.

\smallskip

[Ci] G.~Ciolli. {\it On the quantum cohomology of some Fano threefolds
and a conjecture of Dubrovin.}  Int. Journ. of Math., vol. 16, No. 8 (2005),
823--839. e-arxiv math/0403300.

\smallskip

[Du] B.~Dubrovin. {\it Geometry and analytic theory of Frobenius manifolds.}
Proc. ICM Berlin 1998, vol. 2, 315 -- 326.

\smallskip

[Gi] A.~Givental. {\it A mirror theorem for toric complete
intersections.} Progr. Mat., 160 (1998), 141--175.

\smallskip

[Go1] V.~Golyshev. {\it Riemann--Roch variations.} Izv. Ross. AN, ser. mat.,
65:5 (2001), 3--32.

\smallskip

[Go2] V.~Golyshev. {\it A remark on minimal Fano threefolds.} 
e-arxiv 0803.0031

\smallskip

[He] C.~Hertling. {\it Frobenius spaces and moduli spaces for singularities.}
Cambridge University Press, 2002.

\smallskip

[HeMa] C~Hertling, Yu.~Manin. {\it Weak Frobenius manifolds.}
Int. Math. Res. Notices, 6 (1999), 277--286. e-arxiv
math.QA/9810132.

\smallskip

[Ma1] Yu. Manin. {\it Frobenius manifolds, quantum cohomology, and moduli spaces.}
AMS Colloquium Publications, vol. 47, 1999.

\smallskip

[Ma2] Yu. Manin. {\it Manifolds with multiplication on the tangent sheaf.}
Rendiconti Mat. Appl., Serie VII, vol. 26 (2006), 69--85. e-arxiv math.AG/0502016

\smallskip

[Me1] S.~Merkulov. {\it Operads, deformation theory and $F$--manifolds.}
In: Frobenius manifolds, quantum cohomology, and singularities
(eds. C. Hertling and M. Marcolli), Vieweg 2004,  213--251.
e-arxiv math.AG/0210478.

\smallskip

[Me2] S.~Merkulov. {\it PROP profile of Poisson geometry.} Comm. Math. Phys.,
262:1 (2006), 117--135.
e-arxiv math.AG/0401034.

\smallskip

[Or] D.~Orlov. {\it Derived category of coherent sheaves, and motives.} 
(Russian)  Uspekhi Mat. Nauk  60  (2005),  no. 6(366), 231--232;  translation in  Russian Math. Surveys  60  (2005),  no. 6, 1242--1244 14F05
e-arxiv math/0512620

\smallskip

[Sa]  C.~Sabbah. {\it Isomonodromic deformations and Frobenius manifolds.}
EDP Sciences and Springer, 2007

\smallskip

[Te] C.~Teleman. {\it The structure of 2D semi-simple field theories.} 
e-arxiv 0712.0160

\smallskip

[Z] E. Zaslow. {\it Solitons and helices: the search for a math--physics bridge.}
Comm. Math. Phys. 175, No. 2 (1996), 337--375.

\bigskip
{\it Claus Hertling, Institut f\"ur Mathematik, Universit\"at Mannheim,
A5, 6, 68131 Mannheim, Germany, hertling\@math.uni-mannheim.de

\smallskip

Yuri I. Manin, Northwestern University, Evanston, USA, and 
Max-Planck-Institut f\"ur Mathematik, Bonn, Germany, manin\@mpim-bonn.mpg.de

\smallskip

Constantin Teleman, University of Edinburgh, UK, and UC Berkeley,USA,  

c.teleman\@ed.ac.uk}

\enddocument